\documentclass[a4paper,12pt]{article}
\usepackage{tabularx} 
\usepackage{amsmath,amsthm,verbatim,amssymb,amsfonts,amsopn,amscd}
\usepackage{graphicx} 
\usepackage[margin=1in,a4paper]{geometry} 

\usepackage[utf8]{inputenc}
\usepackage[T1]{fontenc}

\usepackage[UKenglish]{babel}
\usepackage{enumerate}
\usepackage{paralist} 
\usepackage[final]{hyperref} 
\hypersetup{
	colorlinks=true,       
	linkcolor=blue,        
	citecolor=blue,        
	filecolor=magenta,     
	urlcolor=blue         
}

\usepackage{tikz}
\usetikzlibrary{decorations.markings}

\theoremstyle{plain}
\newtheorem{theorem}{Theorem}[section]

\theoremstyle{definition}
\newtheorem{definition}[theorem]{Definition}
\newtheorem{remark}[theorem]{Remark}

\DeclareMathOperator{\calC}{\mathcal{C}}

\DeclareMathOperator{\Dih}{Dih}
\DeclareMathOperator{\HO}{HO}

\newcommand{\IZ}{\mathbb{Z}}
\newcommand{\IQ}{\mathbb{Q}}

\newcommand{\la}{\langle}
\newcommand{\ra}{\rangle}

\title{Towards a geometric theory of integration}
\author{\href{mailto:filip.bar.research@gmail.com}{Filip B\'{a}r}}
\date{8th March 2024}

\begin{document}

\maketitle
\begin{abstract}
  \noindent Integration is the final key step when turning an infinitesimal argument into a result applicable to quantities of finite size. Conceptually, it is about combining infinitesimal contributions to a finite whole. We make a first step towards a geometric theory of integration in the context of Synthetic Differential Geometry (SDG) by analysing the differential aspect of the integration process. Starting from two heuristic principles that combine the idea of differential forms as infinitesimal measures while formalising the process of taking infinitesimal differences at the same time we derive a general notion of differential form as an equivariant map from infinitesimal $n$-cuboids to the base ring coordinatising a line. Besides the familiar differential forms introduced by Cartan we discover two new types. We also discover a new differential operator besides the exterior derivative.  
  Analogous to the relationship between the exterior derivative and the Stokes-Cartan integral theorem, this new operator is linked to the generalised Fundamental Theorem of Calculus in higher dimensions, as discussed in prior research.
  This shows that the Fundamental Theorem is an integral theorem like Stokes-Cartan, but for one of the new types of differential forms.      
\end{abstract}

\section{Introduction} 

\noindent Synthetic Differential Geometry (SDG) is a field that has been born out of Bill's endeavour of providing an axiomatic foundation for continuum physics striving for simplicity and clarity. After decades of work and refinement of the initial ideas (most notably by Anders Kock) we have a naive version of SDG that has proven itself to be excellent at modelling geometric intuition based on arguments using infinitesimals; at least when it comes down to doing algebra with infinitesimal structures and constructions. It has helped us to gain much deeper insights into the nature of infinitesimals, infinitesimal arguments, and infinitesimal constructions in geometry than an approach purely based on intuition would ever permit; very much as Bill has envisaged. (See, for example, the work of Kock, Reyes, Bunge, Wraith, Lavendhomme, Penon, Bell and the author.) However, Bill has also stressed repeatedly that the search for clarity is an ongoing process and he pointed out that an informed modification of the established axioms might be necessary. The aim of this note is to make a case for such a modification.

Despite the successes of SDG, there is one important aspect of infinitesimals that SDG has not been able to capture well, so far: \emph{Integration}. Integration is the final key step when turning an infinitesimal argument, which is of evanescent nature, into a result applicable to quantities of finite size. Conceptually, it is about combining infinitesimal contributions to a finite whole. The infinitesimal contributions are differential forms serving as measures on infinitesimal pieces of space with values in some space of quantities. This space carries the structure of a group; for, as the term \textit{`differential'} suggests, we need to be able to take differences. The integration process itself requires a tiling of a finite piece of space into infinitesimal tiles over which the infinitesimal contributions of the differential form can be composed to a `finite' whole.

The current \emph{integration axiom} in SDG \cite{Kock_Reyes:integration_axiom} does not capture this process at all. All it can provide us with is the existence of the integration result. It does so by using the Fundamental Theorem of Calculus as a retrospective definition of the integral rather than a theorem. In this effective approach to integration the integration process itself remains but an intuition that has no direct formalisation in the current SDG. The aim of this note is to conduct a careful analysis of the differential aspect of the integration process, and use the result of this analysis to make a strong case for the need of a \emph{synthetic theory of integration} that is of \emph{geometric} rather than arithmetic nature. 

For our analysis we shall adopt two heuristic principles at the meta-theory level: 
\begin{enumerate}[(A)]
  \item \label{principle:A} Differential forms are measures on infinitesimal pieces of space. 
  \item \label{principle:B} Differential forms are obtained from taking infinitesimal differences.
\end{enumerate}
Each of these principles will lead us independently to the same new type of differential form that turns out to be conceptually different from the differential forms introduced by Cartan \cite{Cartan:diffforms}. We shall call these new type \emph{Leibniz-Fubini differential forms} for the time being\footnote{This is to honour G.~W.~Leibniz and G.~Fubini. It was Leibniz, who first introduced the concept of a differential in his work on infinitesimal calculus; whereas the defining symmetry of the new type of differential forms follows from Fubini's theorem about the interchangeability of the order of integration.}.
Principle~(\ref{principle:B}) demands the infinitesimal pieces of space our differential forms take as arguments to be infinitesimal quadrilaterals (quads) for 2-forms and, more generally, infinitesimal $n$-cuboids for $n$-forms. (We shall only focus on 2-forms here, as they are sufficient to make our case. The full theory of $n$-forms is currently work in progress and will be presented elsewhere \cite{Bar:extending_Cartan}.) 

The idea to use infinitesimal cuboids as a domain of definition for differential forms is not new in SDG, and has been first proposed in \cite{Kock_et_al:forms_and_integration_in_SDG} (see also   
\cite[ch.~4]{Lavendhomme:SDG}). However, in \cite{Kock_et_al:forms_and_integration_in_SDG} the authors consider parametrised cuboids $D^n\to M$ on a microlinear space $M$, and the set of first-order infinitesimals $D\subset R$ on a ring of line type\footnote{Recall that $R$ is said to be of line type if it satisfies the Kock-Lawvere axiom schema (see, e.g., \cite[ch.~2.1.3]{Lavendhomme:SDG}), which makes precise the intuition that any function restricted to an infinitesimal variety becomes a polynomial.} $R$; whereas we shall define infinitesimal quads using infinitesimally neighbouring points similar to the approach taken in the theory of combinatorial differential forms in \cite[ch.~3]{Kock:Synthetic_Geometry_Manifolds}.

Another consequence of principle~(\ref{principle:B}) is the existence of a differential operator $\delta$ in addition to the exterior derivative $d$. It maps 1-forms to Leibniz-Fubini 2-forms. The key observation is that $\delta^2$ is not zero and maps functions, i.e. 0-forms, to Leibniz-Fubini 2-forms. Assuming the \emph{integration axiom}, i.e. any function $f:[a,b]\to R$ has an antiderivative $F:[a,b]\to R$ such that $F'=f$ and the latter is unique up to a constant, any Leibniz-Fubini 2-form $f(x,y)\,dx\,dy$ on $I=[a_1,b_1]\times[a_2,b_2]$ is of the form $\delta^2 F$ for some function $F:I\to R$. Moreover, this $F$ is a generalised antiderivative of $f$ and can be used to evaluate the integral of $f$ over $I$ as in the generalisation of the Fundamental Theorem of Calculus to 2D \cite{Bar:fundamental_theorem_of_calculus} (see also equation~(\ref{eq:FTC-2D})). In practice we calculate any integral of a 2-form in this way; so any 2-form ends up being identified with a Leibniz-Fubini 2-form during the process of integration, which is unnatural.

Parallel to Stokes-Cartan, which provides an integral formula corresponding to the exterior derivative $d$, we shall show the generalised Fundamental Theorem of Calculus to be the integral theorem corresponding to $\delta$. It comes in a "relative" and an "absolute" version. The relative version is an integral theorem like Stokes-Cartan, but for Leibniz-Fubini forms instead of Cartan forms (see equation~(\ref{eq:relative-FTC-2D})). We expect the absolute version to generalise the combinatorial formula in \cite{Bar:fundamental_theorem_of_calculus} to (certain) manifolds: Given a Leibniz-Fubini 2-form $\alpha$ on a two-dimensional manifold $M$ with four marked vertices $P,Q,R,S$ such that $\alpha = \delta^2 F$ for a function $F:M\to R$, and $M$ admits a tiling into infinitesimal quads compatible with the markings, then
$$
  \int_M \alpha = F(R) - F(S) + F(P) - F(Q) 
$$
With the current integration axiom alone this conjecture cannot be investigated properly, however. Indeed, the only way to integrate 2-forms in SDG as of now is to first pull them back via a parametrisation and then integrate the pullback over a rectangular domain. In this case the above equation becomes a consequence of the generalised Fundamental Theorem of Calculus as proved in \cite{Bar:fundamental_theorem_of_calculus} and the change of variable theorem. To be able to go beyond this case requires a notion of an integral of a 2-form on a manifold in SDG without having to invoke a parametrisation apriori. \emph{Such a notion ought to be provided by a geometric theory of integration.}

Throughout this paper we will be working in the context of naive SDG over a $\IQ$-algebra $R$ that satisfies the following \emph{Kock-Lawvere axioms}: For any natural number $n\geq 1$ let $D(n)\subset R^n$ be the set of first-order infinitesimals 
$$
  D(n) =\{(d_1,\ldots,d_n)\in R^n\mid d_id_j=0,\ 1\leq i,j\leq n\}
$$
Any function $f:D(n)\to R$ is affine, i.e. there are unique $a_k\in R$, $0\leq k\leq n$ such that 
\begin{equation}\label{eq:FTC-2D-manifold}
  f(d_1,\ldots,d_n)=a_0 + \sum_{k=1}^n a_kd_k
\end{equation} 
for all $(d_1,\ldots,d_n)\in D(n)$. For a function $f:R^n\to R$ and $1\leq k\leq n$ we define the \emph{$k$th partial derivative} $\partial_k f(x_0)$ as the uniquely determined $a_k\in R$ when applying the Kock-Lawvere axiom to $d\mapsto f(x_0+d\, e_k)$, where $e_k$ is the $k$-th standard basis vector in $R^n$ and $d\in D=D(1)$.

The $\IQ$-algebra $R$ is assumed to be pre-ordered. Apart from being a pre-ordered ring this means that the pre-order cannot distinguish infinitesimals from $0$: $0\leq d\leq 0$ for all $d\in D$. As a consequence all intervals $[a,b]\subset R$ are considered to be \emph{marked}, i.e. the points $a,b\in R$ are part of the data. Due to $0\leq d\leq 0$ any $f:[a,b]\to R$ can be differentiated on the whole interval including $a$ and $b$. We also assume the integration axiom.

\section{Integration from the infinitesimal viewpoint}\label{sec:integration}

Leibniz's integral notation $\int_a^b f(x)\,dx$ fully reflects the integration process of an infinitesimal argument as outlined in the introduction. For example, when finding the area enclosed by the graph of $f$, the $x$-axis and the vertical lines $x=a$ and $x=b$, the expression $f(x)\,dx$ would be the infinitesimal contribution to that area by a rectangle of height $f(x)$ and infinitesimal width $dx$. The area under the graph of $f$ is obtained by combining all these infinitesimal contributions over a subdivision of $[a,b]$ into infinitely many intervals of infinitesimal width to a finite whole.

In Analysis the differential $dx$ is not taken seriously and attributed a symbolic meaning only; whereas one can interpret $f(x)\,dx$ as a 1-form in the framework of Cartan's formalisation of differential forms in Differential Geometry, or Global Analysis, at least \cite{Cartan:diffforms,Cartan:diffforms_book}. In this formalisation, however, differential forms merely select alternating multilinear maps on the tangent spaces at each point $x$ of the domain we want to integrate over in a continuous way. There is no notion of infinitesimal; $dx$ is but a linear projection map. Integrating such a differential form relies on the integrals from Analysis, where the 1-form $dx$ has but a symbolic meaning again.    

In SDG the infinitesimal expression $f(x)\,dx$ does not only make sense for $dx\in D$, but is also indicative of why the Fundamental Theorem of Calculus holds true. Indeed, the integration axiom guarantees an antiderivative $F$ of $f$, and the Kock-Lawvere axiom for $n=1$ shows 
the differential form $f(x)\,dx$ a differential of $F$:
$$
  F(x+dx) - F(x) = f(x)\,dx
$$
By principle~(\ref{principle:A}) we have that
$$
  \int_{[x,x+dx]}f(u)\,du = f(x)\,dx
$$
Combining all the infinitesimal contributions $f(x)\,dx$ over any subdivision of $[a,b]$ into infinitesimal intervals becomes a telescoping sum that evaluates to $F(b)-F(a)$: 
$$
  \int_{[a,b]}f(x)\,dx = \sum_{x\in[a,b]}\int_{[x,x+dx]}f(u)\,du = F(b) - F(a)
$$
Since we have currently no means to make precise the idea of the infinitesimal subdivision and combination process in SDG, the above argument remains but a heuristic justification of
the \emph{definition of the integral} in SDG via the antiderivative given by the integration axiom. A proper synthetic \emph{geometric theory of integration} should provide the means to make this argument precise.

The crucial observations can only be made when considering two or higher dimensional integrals, however. Consider the integral of a function $f:[a_1,a_2]\times [b_1,b_2] \to R$ over the rectangular domain $[a_1,a_2]\times [b_1,b_2]$ 
$$\int_{[a_1,a_2]\times [b_1,b_2]} f(x,y)\,dx\,dy$$
The integrand is a differential form $f(x,y)\,dx\,dy$ that \emph{cannot be a 2-form in the sense of Cartan}. Indeed, we have $f(x,y)\,dx\wedge dy = - f(x,y)\,dy\wedge dx$, but $f(x,y)\,dx\,dy$ is clearly symmetric in $dx$ and $dy$. 

As in the 1D case we can attempt to write the integrand as a differential of a function $F:[a_1,a_2]\times [b_1,b_2] \to R$. Let $dx,dy\in D$. Applying Kock-Lawvere for $n=1$ to the first and second argument in succession yields
$$
	F(x+dx,y+dy) = F(x,y) + \partial_1 F(x,y) dx + \partial_2 F(x,y) dy + \partial_2\partial_1 F(x,y)\,dx\,dy
$$   
We read off that $F$ has to satisfy the PDE $\partial_2\partial_1 F(x,y) = f(x,y)$. Applying the integration axiom first to the $y$ and then to the $x$-variable yields that such a function $F$ always exists. In a next step we can make $f(x,y)\,dx\,dy$ the subject and replace the partial derivatives with their respective differentials:
\begin{equation}\label{eq:Leibniz-Fubini-2-form}
  	f(x,y)\,dx\,dy = F(x+dx,y+dy) - F(x,y+dy) + F(x,y) - F(x+dx,y)
\end{equation}
Applying principle~(\ref{principle:A}) we find
\begin{align*}
    \int_{[x,x+dx]\times [y,y+dy]} f(u,v)\,du\,dv &= f(x,y)\,dx\,dy \\
    &=  F(x+dx,y+dy) - F(x,y+dy) + F(x,y) - F(x+dx,y)
\end{align*}
The signs at each vertext can be visualised as follows
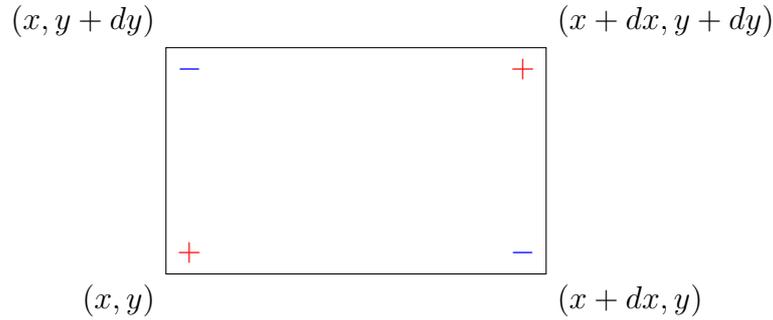
\begin{center}
  \begin{tikzpicture}
    \draw (0,0) rectangle (5,3);
    
    \node at (0,0) [below left] {$(x,y)$};
    \node at (5,0) [below right] {$(x+dx,y)$};
    \node at (5,3) [above right] {$(x+dx,y+dy)$};
    \node at (0,3) [above left] {$(x,y+dy)$};
    
    \node at (0,0) [above right, red] {$+$};
    \node at (5,0) [above left, blue] {$-$};
    \node at (5,3) [below left, red] {$+$};
    \node at (0,3) [below right, blue] {$-$};
  \end{tikzpicture}    
\end{center}
Moreover, the above combinatorial formula is composable:
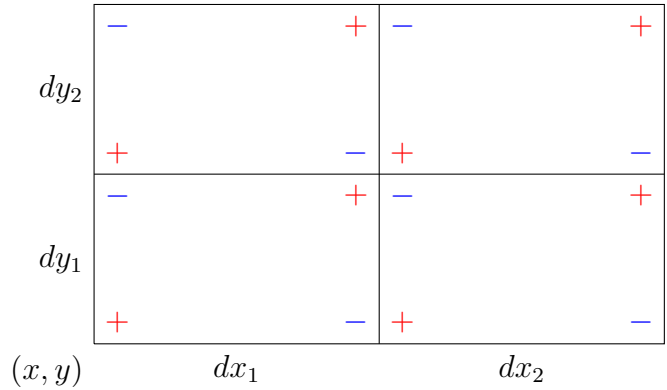
\begin{center}
  \begin{tikzpicture}[scale=1.5]
    \draw (0,0) -- node[midway, below] {$dx_1$} (2.5,0) -- node[midway, below] {$dx_2$} (5,0) -- (5,3) -- (0,3) -- node[midway, left] {$dy_2$} (0,1.5) -- node[midway, left] {$dy_1$} (0,0) node[below left] {$(x,y)$};
    
    \draw (0,1.5) -- (5,1.5);
    \draw (2.5,0) -- (2.5,3);
    
    \node at (0,0) [above right, red] {$+$};
    \node at (5,0) [above left, blue] {$-$};
    \node at (5,3) [below left, red] {$+$};
    \node at (0,3) [below right, blue] {$-$};
    
    \node at (2.5,0) [above left, blue] {$-$};
    \node at (2.5,0) [above right, red] {$+$};
    
    \node at (2.5,1.5) [below right, blue] {$-$};
    \node at (2.5,1.5) [below left, red] {$+$};
    \node at (2.5,1.5) [above left, blue] {$-$};
    \node at (2.5,1.5) [above right, red] {$+$};
    
    \node at (2.5,3) [below right, blue] {$-$};
    \node at (2.5,3) [below left, red] {$+$};
    
    \node at (0,1.5) [below right, blue] {$-$};
    \node at (0,1.5) [above right, red] {$+$};
    
    \node at (5,1.5) [above left, blue] {$-$};
    \node at (5,1.5) [below left, red] {$+$};
    
  \end{tikzpicture}  
\end{center}

Assuming a tiling of $[a_1,a_2]\times [b_1,b_2]$ into infinitesimal rectangles, summing the infinitesimal contributions over the hypothetical partition yields a telescoping sum once again. We recover the combinatorial formula of the generalised Fundamental theorem of Calculus in 2D \cite{Bar:fundamental_theorem_of_calculus} that only involves the outermost vertices:
\begin{equation}\label{eq:FTC-2D}
  	\int_{[a_1,a_2]\times [b_1,b_2]} f(x,y)\,dx\,dy =  F(a_2,b_2) - F(a_1,b_2) +F(a_1,b_1) - F(a_2,b_1)
\end{equation}
We summarise the three main observations made in this section:
\begin{enumerate}
  \item A geometric theory of integration should provide means to discuss tilings into infinitesimal shapes like intervals and rectangles within SDG, as well as means how to compose infinitesimal contributions over such a tiling. Together with the integration axiom we can then derive the generalised Fundamental Theorem of Calculus, rather than using it to effectively define the integral via the integration axiom, retrospectively.
  \item The integrand $f(x,y)\,dx\,dy$ is not a 2-form in the sense of Cartan, but an example of a \emph{Leibniz-Fubini 2-form}. 
  \item Equation~(\ref{eq:Leibniz-Fubini-2-form}) shows that there is an operation mapping a function to a Leibniz-Fubini 2-form. Moreover, the Fundamental Theorem of Calculus in 2D is equivalent to stating that every Leibniz-Fubini 2-form on a rectangular domain is of this form.
\end{enumerate}

Tackling observation~(1) is the main challenge in developing a synthetic geometric theory of integration. It is current work in progress and beyond the scope of this note. We shall focus on exploring principle~(\ref{principle:B}) to define a general notion of differential form that encompasess both Cartan 2-forms and Leibniz-Fubini 2-forms. We will also define the symmetric differential operator $\delta$ and see that the operation mentioned in the third point is indeed $\delta^2$.

\section{Three types of differential 2-forms}\label{sec:3_types_of_2-forms}

To define a general notion of a 1-form on $R^n$ we can allude to principle~(\ref{principle:B}) and Leibniz's example of a 1-form $f(x)\,dx$ on $R$ as a differential of a function $F(x)$. Since we want to define differential forms on points, we require pairs of points that are infinitesimally close to each other as the domain of definition:
$$
  R^n\la 2\ra = \{(P,Q)\in (R^n)^2\mid Q-P\in D(n)\} 
$$
Adopting the terminology from algebraic geometry this space is referred to as the first neighbourhood of the diagonal in SDG literature (see e.g. \cite{Kock:Synthetic_Geometry_Manifolds}). Our chosen notation aligns with that of infinitesimal structures \cite{Bar:second_order_affine_structures,Bar:2nd-order_i-groups}, as this is indicative of possible further generalisations of this notion of differential form to other contexts.

Using $R^n\la 2\ra\cong R^n\times D(n)$ and Kock-Lawvere we see that any map $\alpha:R^n\la 2\ra\to R$ is of the form 
$$
  \alpha(P,Q) = a_0(P) + \sum_{k=1}^n a_k(P)(Q_k-P_k)
$$
respectively, using the dot product on $R^n$
$$
  \alpha(P,Q) = a_0(P) + X(P)\bullet (Q-P)
$$
for a unique vector field $X:R^n\to R^n$ and function $a_0:R^n\to R$. By principle~(\ref{principle:B}) we want 1-forms to formalise taking differences, so we demand $\alpha(P,P)=0$ for all $P\in R^n$, and hence $a_0=0$. With this we recover the notion of a \emph{combinatorial 1-form} in SDG \cite[ch.~3]{Kock:Synthetic_Geometry_Manifolds} (and their local representation in a chart). In particular, any 1-form $\alpha$ on $[a,b]\subset R$ is of the form $f(x)\,dx$ for a uniquely determined function $f:[a,b]\to R$. 

Taking the differential of a function yields an operation $\delta$ that maps functions to 1-forms. It is defined by
$$
  \delta f(P,Q) = f(Q) - f(P), \quad P,Q\in R^n\la 2\ra
$$
and the corresponding vector field is the gradient $\nabla f$ of $f$. 

According to principle~(\ref{principle:B}) we ought to construct a 2-form out of a 1-form by taking its differential. Firstly, we need to define what it means for two pairs $(P,Q), (R,S)\in R^n\la 2\ra$ to be infinitesimal neighbours. The simplest way is to do this componentwise, i.e. to require $(P,R), (Q,S)\in R^n\la 2\ra$. Geometrically, this yields an infinitesimal quad $PQSR$ (labeled anti-clockwise).
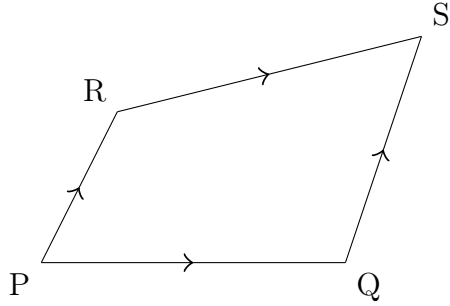
\begin{center}
  \begin{tikzpicture}
    \begin{scope}[decoration={
      markings,
      mark=at position 0.5 with {\arrow[scale=2]{to}}}
      ] 
      \draw[postaction={decorate}] (0,0) node[below left] {P} -- (4,0) node[below right] {Q};
    
      \draw[postaction={decorate}] (1,2) node[above left] {R} -- (5,3) node[above right] {S};
    
      \draw[postaction={decorate}] (4,0) -- (5,3);
    
      \draw[postaction={decorate}] (0,0) -- (1,2);   
    \end{scope}
  \end{tikzpicture}    
\end{center}

We have two means of taking the differential of a 1-form $\alpha$ on such a quad: vertically
$$
  \delta_v\alpha(P,Q,S,R) = \alpha(R,S) - \alpha(P,Q)
$$
or horizontally
$$
  \delta_h\alpha(P,Q,S,R) = \alpha(Q,S) - \alpha(P,R)
$$
Since there is no geometric reason to prefer one over the other, we take the symmetrisation 
\begin{equation}\label{eq:symmetric_diff}
    \delta\alpha(P,Q,S,R) = \frac{1}{2}(\delta_h\alpha(P,Q,S,R) + \delta_v\alpha(P,Q,S,R))
\end{equation}
and antisimmetrisation of both differentials
\begin{equation}\label{eq:exterior_diff}
    d\alpha(P,Q,S,R) = \delta_h\alpha(P,Q,S,R) - \delta_v\alpha(P,Q,S,R)
\end{equation}
We have dropped the factor $1/2$ in the definition of the \emph{boundary differential} $d\alpha$, for we recover the \emph{exterior derivative} of a 1-form in the sense of Cartan this way. Indeed, by applying principle~(\ref{principle:A}) we make the same observation as in \cite{Kock_et_al:forms_and_integration_in_SDG}, namely that the exterior derivative is just the infinitesimal version of the Stokes-Cartan theorem.    
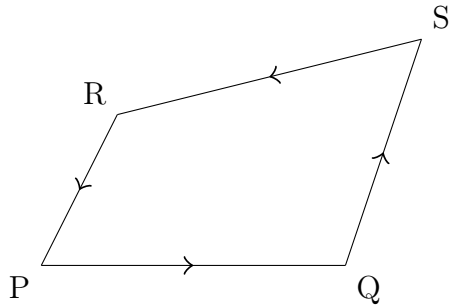
\begin{center}
  \begin{tikzpicture}
    \begin{scope}[decoration={
      markings,
      mark=at position 0.5 with {\arrow[scale=2]{to}}}
      ] 
      \draw[postaction={decorate}] (0,0) node[below left] {P} -- (4,0) node[below right] {Q};
    
      \draw[postaction={decorate}] (5,3) node[above right] {S} -- (1,2) node[above left] {R};
    
      \draw[postaction={decorate}] (4,0) -- (5,3);
    
      \draw[postaction={decorate}] (1,2) -- (0,0);   
    \end{scope}
  \end{tikzpicture}    
\end{center}

\begin{equation}\label{eq:Stokes-Cartan-infinitesimal}
  \begin{split}
    \int_{PQSR} d\alpha &= d\alpha(P,Q,S,R) \\
         &= \alpha(P,Q) + \alpha(Q,S) + \alpha(S,R) + \alpha(R,P) \\ 
      &=\int_{\partial PQSR} \alpha
  \end{split} 
\end{equation}
This justifies calling the operator $d$ the boundary differential.

The operator $\delta$ is the \emph{symmetric differential operator}. Indeed, for $\alpha = \delta F$ we recover the combinatorial formula~(\ref{eq:Leibniz-Fubini-2-form})
$$
  \delta\alpha(P,Q,S,R) = \delta^2 F(P,Q,S,R) = F(S) - F(R) + F(P) - F(Q)
$$
Applying principle~(\ref{principle:A}) this yields the infinitesimal geometric version of the Fundamental Theorem of Calculus for Leibniz-Fubini 2-forms (cf. equation~(\ref{eq:FTC-2D})).

\begin{equation}\label{eq:FTC-geometric-infinitesimal}
 \int_{PQSR} \delta^2 F = \delta^2 F(P,Q,S,R) = F(S) - F(R) + F(P) - F(Q) 
\end{equation}

Defining Cartan and Leibniz-Fubini differential forms in general can be done via their behaviour under the action of the symmetry group $\Dih_4$ of a square by re-labeling of the infinitesimal quads. Recall that $\Dih_4$ has two generators $r$ and $s$, where the intended interpretation of $r$ is a rotation by $\pi/2$ and the intended interpretation of $s$ is a reflection. We define the action by having $r$ act on the labels as a cyclic permutation to the right, and $s$ swapping the labels of the bottom right and top left vertex; i.e. as a reflection of the square along the diagonal through the bottom left and top right vertex 
$$
  r_*(PQRS) = SPQR, \qquad s_*(PQRS) = PSRQ
$$
There are three non-trivial group homomorphisms $\sigma_j: \Dih_4\to \{-1,1\}$ determined as follows:
\begin{align*}
  \sigma_{C}(r) &= 1, & \sigma_{C}(s) &= -1 \\
  \sigma_{LF}(r) &= -1, & \sigma_{LF}(s) &= 1 \\
  \sigma_{N}(r) &= -1, & \sigma_{N}(s) &= -1
\end{align*} 

Before we can give the definition of a general differential 2-form we need to review what we mean by an \emph{infinitesimal quad}. So far, we have treated a quad as a purely combinatorially defined notion based on the first-order neighbourhood relation. However, we have no guarantee that the infinitesimal quads are planar. In fact, it is not hard to construct examples of maps that map planar infinitesimal quads to non-planar ones.

We have to accept non-planar quads, but we do require that the four vertices of the infinitesimal quad enclose a \emph{volume of zero}. This can be justified by fluid mechanics, or more generaly, by classical field theory. In 3D a Cartan 2-form is measuring an infinitesimal current, like the rate of flow of mass, or rate of flow of charge through an infinitesimal quad, for example. If the quad encloses a volume that is not provably equal to zero, then it could enclose sources or sinks of the current. In this case there is a net infinitesimal contribution to the current through the surface by the surface itself, which would confuse volume with surface integrals and make Stokes' theorem fail in 3D. 

\begin{definition}[Infinitesimal quad]\label{def:infinitesimal_quad} 
  A $4$-tuple of points $PQRS$ forms an \textbf{\emph{infinitesimal quad}} on $R^n$ if
  \begin{enumerate}
    \item[(Q1)] $(P,Q), (Q,R), (R,S), (S,P)\in R^n\la 2\ra$
    \item[(Q2)] $b[Q-P,S-P,R-S]=0$ for any alternating trilinear form on $R^n$
  \end{enumerate}
  The set of all infinitesimal quads on $R^n$ is denoted by $\calC_2(R^n)$. 
\end{definition}

Infinitesimal quads are thus combinatorial quads that are "thin" or "measurably flat".

\begin{definition}[Differential 2-forms]\label{def:2-forms} 
  Let $\sigma:\Dih_4\to \{-1,1\}$ be a non-trivial group homomorphism. A \textbf{\emph{differential 2-form}} on $R^n$ is a map $\omega:\calC_2(R^n)\to R$ such that
  $$
    \omega(t_*PQRS) = \sigma(t)\,\omega(PQRS)
  $$
  for all $t\in \Dih_4$ and $PQRS\in \calC_2(R^n)$. 
\end{definition}

We call a differential 2-form $\alpha$ \emph{Cartan} if $\sigma = \sigma_C$, \emph{Leibniz-Fubini} if $\sigma = \sigma_{LF}$, and \emph{Nieuwentijdt}\footnote{This is to honour B. Nieuwentijdt, a contemporary of Leibniz who also attempted an axiomatisation of infinitesimals. In contrast to Leibniz Nieuwentijdt argued that infinitesimals should square to zero \cite{Bell:continuous_and_infinitesimal}.} if $\sigma = \sigma_N$. It is not hard to see that $d\alpha$ is indeed a Cartan 2-form, and $\delta^2 F$ a Leibniz-Fubini 2-form according to the above definitions. Nieuwentijdt forms have only been discovered by noting that $\Dih_4$ (and more generally, all the hyperoctahedral groups) have three subgroups of index $2$. However, we are still looking for examples where such forms arise in applications. We shall not discuss them here. 

\begin{remark}\label{rem:1-forms}
  Note that differential 1-forms can be defined via the above group action approach as well. Firstly, any alternating bilinear form on $R^n$ will vanish on $Q-P$ for $(P,Q)\in R^n\la 2\ra$. We thus have $\calC_1(R^n) = R^n\la 2\ra$. Secondly, the symmetry group of a line segment is $\{-1,1\}$, so the only non-trivial group homomorphism is the identity map; but this is just saying that a 1-form is alternating, which is equivalent to saying that it vanishes on the diagonal for a ring with a characteristic zero.   
\end{remark}

We conclude this section by giving representations of Cartan and Leibniz-Fubini 2-forms on $R^n$. They can be obtained by applying the Kock-Lawvere axiom, the symmetry conditions and the zero-volume condition of infinitesimal quads. We shall give the representations for general infinitesimal quads. By restricting to infinitesimal parallelograms, or rectangles, the representations simplify; particularly for the Leibniz-Fubini forms.

Cartan 2-forms have the simplest representation. Any such 2-form $\omega$ is equivalent to a map $b$ from $R^n$ into the space of alternating bilinear forms on $R^n$
$$
  \omega(PQRS) = b_P[Q-S.R-P]
$$
The representation of Leibniz-Fubini 2-forms is more involved. Any such $\omega$ determines a map $A$ from $R^n$ to the space of linear forms as well as a map $a$ from $R^n$ to the space of bilinear forms
\begin{align*}
    \omega(PQRS) = &A_P[R-S+P-Q] \\
        &+\frac{1}{2}\bigl(\partial A_P[Q-P,R-Q]+\partial A_P[S-P,R-S]\bigr)\\
        &+ a_P[R-P]^2-a_P[Q-S]^2	
\end{align*}
Here $\partial A_P$ denotes the derivative of $A$ at $P$, which is a bilinear form. Conversely, any two such maps determine a Leibniz-Fubini 2-form using the mapping rule above. On an infinitesimal parallelogram $PQRS$ we have $R-S=Q-P$ and $R-Q=S-P$. In this case the representation simplifies to
$$
  \omega(PQRS) = b_p[Q-P,S-P],
$$
where $b$ is a map from $R^n$ to the space of symmetric bilinear forms on $R^n$. In particular, we have 
$$
  \delta^2 F(PQRS) = \partial^2F(P)[Q-P,S-P]
$$
on parallelograms.

Although, we have stated all definitions for 2-forms for $R^n$ only, this has a straight forward generalisation to 2-forms on formally open subsets and hence to 2-forms on a manifold. To define $n$-forms the definition of infinitesimal quads needs to be extended to infinitesimal $n$-cuboids, and $\Dih_4$ needs to be replaced by the hyperoctahedral group $\HO(n)$. (Note that $\Dih_4=\HO(2)$) As the latter has exactly three index-two subgroups for all $n\geq 2$, $n$-forms also split into the three types introduced here.   


\section{Towards a geometric theory of integration}\label{sec:integration_2-forms}

In current SDG any integration has to be done over a cartesian product of intervals, i.e. a rectangular axis-parallel cuboid in $R^n$. This is due to the effective definition of the integral via antiderivatives and the absence of an independent integral notion. Although this is sufficient for all practical purposes of computing integrals, it is lacking when one wishes to conduct deeper conceptual investigations. In particular, differential forms always need to be pulled back via parametrisations before they can be integrated. 

Any map between manifolds preserves the first-order neighbourhood relation \cite[ch.~2.1]{Kock:Synthetic_Geometry_Manifolds}. Moreover, Kock-Lawvere guarantees that first-order differentials are always mapped linearly. Any mapping between manifolds will thus map infinitesimal quads to infinitesimal quads, and the pullback operation by precomposition maps differential 2-forms to differential 2-forms. It is easy to see that the type of the 2-form is preserved under pullback.

Our prior analysis of the differential aspect of integration now reveals that \emph{the practice of pulling back differential forms and then integrating is geometrically flawed in the current theoretical framework.} Indeed, only Leibniz-Fubini 2-forms can be integrated. Integrating a Cartan 2-form means to pull it back to a rectangular domain in $R^2$, formally replace it with a Leibniz-Fubini 2-form, and then integrate. (Note that the only 2-form that is Cartan and Leibniz-Fubini at the same time is the zero map.)

On the other hand, the infinitesimal geometric versions of \emph{Stokes-Cartan (\ref{eq:Stokes-Cartan-infinitesimal}) and the Fundamental Theorem of Calculus (\ref{eq:FTC-geometric-infinitesimal}) are both composable.} Given a 2D (marked) manifold $M$ with a boundary $\partial M$ (and possibly corners) that admits a tiling into infinitesimal quads, the infinitesimal versions combine to the integral theorems on the whole manifold provided there is some finiteness condition like that $M$ is compact. In the case of the theorem of Stokes, this is how it is usually "proven" when taking a physics or engineering approach to the subject. 

In fact, we can formulate a "relative" version of the Fundamental theorem as an integral theorem for only one application of the symmetric differential $\delta$. Let $\alpha$ be a 1-form. At the infinitesimal level we find

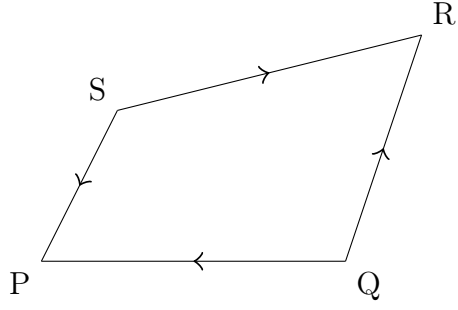
\begin{center}
  \begin{tikzpicture}
    \begin{scope}[decoration={
      markings,
      mark=at position 0.5 with {\arrow[scale=2]{to}}}
      ] 
      \draw[postaction={decorate}] (4,0) node[below right] {Q} -- (0,0) node[below left] {P};
    
      \draw[postaction={decorate}] (1,2) node[above left] {S} -- (5,3) node[above right] {R};
    
      \draw[postaction={decorate}] (4,0) -- (5,3);
    
      \draw[postaction={decorate}] (1,2) -- (0,0);   
    \end{scope}
  \end{tikzpicture}    
\end{center}

\begin{equation}\label{eq:relative-FTC-2D-infinitesimal}
    \begin{split}
      \int_{PQRS} \delta \alpha &= \delta \alpha(PQRS) \\
      &= \frac{1}{2}\bigl(\alpha(Q,P) + \alpha(Q,R) + \alpha(S,R) + \alpha(S,P)\bigr) \\
      &= \frac{1}{2}\Bigl(\int_{QP}\alpha + \int_{QR}\alpha + \int_{SR}\alpha + \int_{SP}\alpha \Bigr)
    \end{split}
\end{equation}
Note that this formula is also composable. It can be extended to 2D manifolds $M$ with boundary $\partial M$ that admit a tiling into infinitesimal quads. However, the boundary of $M$ has to be marked by four vertices compatible with the infinitesimal tiling. If we denote the four vertices by $PQRS$ again, then the corresponding integral theorem reads
\begin{equation}\label{eq:relative-FTC-2D}
  \int_{M} \delta \alpha = \frac{1}{2}\Bigl(\int_{QP}\alpha + \int_{QR}\alpha + \int_{SR}\alpha + \int_{SP}\alpha \Bigr)
\end{equation}
In contrast to Stokes-Cartan the relative Fundamental Theorem can be iterated leading to the "absolute" Fundamental Theorem as stated in the introduction (\ref{eq:FTC-2D-manifold}). This is because $\delta^2$ is not the zero map. 

For an $n$-form we require $2^n$ vertices to be marked, i.e. the $n$-dimensional manifold being a smooth representation of the $n$-cube graph. The integrals on the right-hand side range over the $2n$ facets formed from the vertices in correspondence with the $(n-1)$-cube  facets of the $n$-cube. The "absolute" version of the Fundamental Theorem would then be the generalisation of the Fundamental Theorem of Calculus in $n$ dimensions as stated in \cite{Bar:fundamental_theorem_of_calculus} to manifolds.  

Currently it is not clear whether the Fundamental Theorem generalises to manifolds $M$ or more general spaces beyond those which admit a parametrisation with a domain given by a product of intervals (in which case it holds true due to the change of variable theorem); but the conceptual analysis done in this note hopefully suggests that it is a worthwhile enterprise to go and find out.

\section{Conclusion and Outlook}\label{sec:Conclusion}

Starting from two heuristic principles that combine the idea of differential forms as infinitesimal measures while formalising the process of taking infinitesimal differences at the same time we have derived a general notion of differential form as an equivariant map from infinitesimal $n$-cuboids to $R$. The hyperoctahedral group acts on the $n$-cuboids by re-labeling, whereas $R$ comes equipped with a $\IZ_2$-action by point reflection at $0$. For $n\geq 2$ there are three possible types of differential $n$-forms, which we have named Cartan, Leibniz-Fubini, and Nieuwentijdt for the time being\footnote{This naming convention is to be understood as temporary. Cartan forms could be very well called flux forms as they measure fluxes (related to transport phenomena) in their respective dimension. Leibniz-Fubini are strongly related to measuring content. However, since the geometric meaning of the Leibniz-Fubini and Nieuwentijdt forms is not entirely clear to the author yet, we adopted the current naming convention free of a geometric connotation.}. 

The Cartan forms are well-known from classical differential geometry, albeit in an algebraic rather than the geometric form presented here. A Cartan 2-form induces a simplicial 2-form, and these are known to be equivalent to classical differential 2-forms in the sense of Cartan \cite[thm.~4.7.1]{Kock:Synthetic_Geometry_Manifolds}. Extending classical differential forms to Cartan forms requires an affine connection on points \cite[ch.~2.3]{Kock:Synthetic_Geometry_Manifolds}. Still, the local representation of Cartan 2-forms shows that there is a bijection between Cartan 2-forms and classical differential 2-forms on manifolds.

The Leibniz-Fubini forms are new. They are the differential forms that can be integrated directly and appear to measure content. They have been present since Leibniz (at least), but seem to have been overlooked so far. Unlike Cartan forms it is currently not clear what algebraic construction in classic Differential Geometry corresponds to Leibniz-Fubini forms; in particular, as the 2-forms contain both linear and bilinear terms when evaluated on infinitesimal quads, but are symmetric bilinear on parallelograms. The Nieuwentijdt forms are also new, but no interesting examples or applications have been found yet.

We have found two natural differential operators on differential forms, the boundary differential $d$ and the symmetric differential $\delta$. The former is the infinitesimal version of the exterior derivative mapping 1-forms to Cartan 2-forms, while the latter is new and maps 1-forms to Leibniz-Fubini 2-forms. They both agree on 0-forms, i.e. functions. Unlike $d$ the differential $\delta^2$ is not zero.

To each differential operator corresponds an integration theorem. The boundary differential is the infinitesimal version of the well-known Stokes-Cartan theorem, and the symmetric differential operator turns out to be the infinitesimal version of the generalised relative Fundamental Theorem of Calculus, which is new. Unlike Stokes-Cartan the relative Fundamental Theorem can be iterated resulting in the Fundamental Theorem of Calculus in \cite{Bar:fundamental_theorem_of_calculus} generalised to certain manifolds. Whether there is a differential operator corresponding to Nieuwentijdt forms (possibly with an accompanying integration theorem) is currently unknown.

Besides the Nieuwentijdt forms, which are still a mystery, there are two main open problems as a consequence of the present work. The first open problem is related to the title of this paper: \emph{develop a theory of tilings of spaces into infinitesimal $n$-cuboids with an accompanying notion of composition of differential $n$-forms over such tilings.} This is the hard problem that needs solving to be able to establish a synthetic geometric theory of integration within SDG.

The second open problem is to \emph{ prove the relative Fundamental Theorem of Calculus on manifolds, or other suitable spaces.} The change of variable theorem yields a proof for spaces which admit a parametrisation by a cartesian product of intervals in $R^n$. In this case the class of parametrisations defines the geometry of that space and is part of the structure. To be able to go beyond this case within SDG requires developing an answer to problem one. Tackling this problem in classic Differential Geometry requires constructing a corresponding theory of Leibniz-Fubini forms. The author is grateful for any pointers to differential geometry literature or applications, where such or similar integral theorems could have been proven in hindsight.    


\end{document}